\documentclass[15pt]{article}
 \usepackage{amsfonts,amssymb,amsmath,graphicx}
 \usepackage[cp866]{inputenc}
 \usepackage[english]{babel}
 \usepackage[active]{srcltx}

\newtheorem{theorem}{Theorem}

\author{{\LARGE A.I. Novikova}\thanks{Partially supported by grant RFFR 11-01-00614-a}}
\title{Lyapunov theorem for $q$-concave Banach spaces}
\date{\today}

\begin{document}
\Large
\maketitle

{\large {\bf Abstract.} Generalization of Lyapunov convexity theorem is proved for vector measure with values in Banach spaces with unconditional bases, which are $q$-concave for some $q<\infty.$}


{\bf Introduction.}
Let $X$ be a Banach space, $(\Omega,\Sigma)$ be a measure space, where $\Omega$ is a set and $\Sigma$ is a $\sigma$-algebra of subsets of $\Omega.$ If $m:\Sigma\rightarrow X$ is a $\sigma$-additive $X$-valued measure, then the range of $m$ is the set $m(\Sigma)=\{m(A): \ A\in\Sigma\}.$

The measure $m$ is {\it non-atomic} if for every set $A\in\Sigma$ with $m(A)>0,$ there exist $B\subset A,B\in\Sigma$ such that $m(B)\neq0$ and $m(A \backslash B)\neq0.$

According to the famous Lyapunov theorem ~\cite{Lyap} the range of $\mathbb{R}^n$-valued non-atomic measure $\mu$ is convex. However, this theorem can't be generalized directly for infinite dimensional case: for every infinite dimensional Banach space $X$ there is $X$-valued non-atomic measure (of bounded variation) whose range is not convex ~\cite[Corollary 6]{DU}.

We will call $X$-valued measure {\it a Lyapunov measure} if the closure of its range is convex. And Banach space $X$ is {\it Lyapunov space} if every $X$-valued non-atomic measure is Lyapunov.

Although direct generalization can not be done the following result was obtained in ~\cite{U}.

\begin{theorem} (Uhl). Let $X$ have Radon-Nikodym property. Then any $X$-valued measure of bounded variation is Lyapunov.
\end{theorem}

For measures with unbounded variation this theorem is no longer true. There exist non-atomic measures with values in Hilbert space which are not Lyapunov. It follows that also spaces containing isomorphic copies of $l_2$ are not Lyapunov. In particular, all $L_p[0,1]; \ C[0,1]; \ l_\infty.$ Never the less it was proved ~\cite{lp} that the sequence spaces $c_0$ and $l_p,1\leq p<\infty,p\neq2$ are examples of Lyapunov spaces.

On the other hand, using compactness argument, it was proven in ~\cite{SS} that in a Banach space with unconditional basis any non-negative (with respect to the order induced by the basis) non-atomic measure is Lyapunov.

Recall that for Banach spaces $X,Y$ and $Z$, operator $T:X\rightarrow Y$ is said to be {\it $Z$-strictly singular} if it is not an isomorphism when restricted to any isomorphic copy of $Z$ in $X.$

We say that a linear operator $T\ :\ L_p(\mu) \rightarrow X$ is {\it narrow} if for every $\epsilon > 0$ and every measurable set $A \subset[0, 1]$ there exists $x \in L_p$ with
$x^2 = \mathbb{I}_A$, and $\int_{[0,1]}x d\mu = 0$, so that $\|Tx\| < \epsilon$ (we call such an $x$ a {\it mean zero sign}).

In ~\cite[Theorem B]{MPRS} it was shown that for every $1<p<\infty$ and every Banach space $X$ with an unconditional basis, every $l_2$-strictly singular operator $T:L_p \rightarrow X$ is narrow.

{\bf Main Result.}
The following theorem is a generalization of the result from ~\cite{lp}, presented above.

\begin{theorem}. Let X be a q-concave ($q<\infty$) Banach space with an unconditional basis, which doesn't contain isomorphic copy of $l_2.$
Then X is a Lyapunov space.
\end{theorem}

{\it Proof.} Assume the contrary: X is not Lyapunov, i.e. there is a non-atomic measure $\mu$ with values in X such that the closure of its range is not convex. Then by lemma 3 ~\cite{lp} there exists $(\Omega,\Sigma,\lambda),$ nonnegative measure $\lambda:\Sigma\rightarrow \mathbb{R}: \ \forall A\in\Sigma \ 0\leq\lambda(A)\leq Const\|\mu(A)\|;$ and bounded operator $T:L_\infty(\Omega,\Sigma,\lambda)\rightarrow X$ so that
$$
\begin{array}{l}
\cdot \; T:(L_\infty,w^*)\rightarrow(X,w) \ \mbox{is continous}, \ T(\mathbb{I}_A)=\mu(A);\\
\cdot \; \exists \ \epsilon>0 \ \forall \ \mbox{mean zero sign} \ f\in L_\infty \; \|Tf\|\geq\delta\lambda(supp f).
\end{array}
$$

By factorization theorem ~\cite[1.d.12]{LT2} T can be factorized through $L_q$ (as T is q-concave), i.e. $T=ST_1,$ where
$T_1:L_{\infty}(\Omega,\lambda)\rightarrow L_q(\Omega,\nu)$ is the formal identity map, positive and continuous; and $S:L_q(\Omega,\nu)\rightarrow X$ is bounded. Then
$$
\|\mu(A)\|_X=\|S(\mathbb{I}_A)\|_X\leq\|S\|\cdot\|\mathbb{I}_A\|_{L_q(\nu)}=\|S\|\cdot\nu^{1/q}(A).
$$
Thus there exist $C>0$ so that for all $A\subset\Omega$ we have $0\leq\lambda(A)\leq C\nu^{1/q}(A).$
Then by the Radon-Nikodym theorem we have
$$
\lambda(A)=\int\limits_\Omega y(t)\mathbb{I}_A d\nu,
$$
where $y\in L_1(\Omega,\nu)$ is positive a.e.

Let us choose $\Omega_0\subset\Omega$ of positive measure $\nu$ so that $0<a\leq y(t)\leq b<\infty, \ t\in\Omega_0,$ for some $a,b.$ Fix $\epsilon>0$ and find $y_\epsilon(t)=\sum\limits_{i=1}^n b_i\mathbb{I}_{B_i},$ where ${B_i}$ is partition of $\Omega_0$ and $b_i\geq0,$ so that $\|y(t)-y_\epsilon(t)\|_{2,\nu,\Omega_0}\leq\epsilon.$
Choose some $B_i$ with $\nu(B_i)>0$ and consider arbitrary subset $C\subset B_i, \ C\in\Sigma.$ Then
$$
\begin{array}{c}
|\lambda(C)-b_i\nu(C)|=|\int y(t)\mathbb{I}_C(t)d\nu-\int b_i\mathbb{I}_C(t)d\nu|=\vspace{0.3cm}\\
=|\int(y(t)-y_\epsilon(t))\mathbb{I}_C(t)d\nu|\leq \vspace{0.3cm}\\
\leq \|y-y_\epsilon\|_{2,\lambda}\nu(B_i)^{1/2}\leq\epsilon\nu(B_i)^{1/2}.
\end{array}
$$

Then for any arbitrary mean zero (w.r.t. $\nu$) sign $x$ on $B_i\subset\Omega_0$ we have:
$$
\begin{array}{l}
0=b_i\int\limits_{B_i} x(t)d\nu=b_i(\nu(x=1)-\nu(x=-1))\geq  \vspace{0.2cm}\\
\hspace{2.3cm} \geq \lambda(x=1)-\lambda(x=-1)-2\epsilon\nu(B_i)^{1/2}; \vspace{0.3cm}\\
0=b_i\int\limits_{B_i} x(t)d\nu=b_i(\nu(x=1)-\nu(x=-1))\leq  \vspace{0.2cm}\\
\hspace{2.3cm} \leq \lambda(x=1)-\lambda(x=-1)+2\epsilon\nu(B_i)^{1/2};
\end{array}
$$
i.e.
$$
|\int\limits_{B_i} x(t)d\lambda|\leq{2\epsilon}\lambda(B_i)^{1/2}.
$$
Thus we can find a mean zero (w.r.t. $\lambda$) sign $x_\epsilon,$ $supp(x_\epsilon)\subset supp(x)=B_i$ such that
$$b_i\nu(x_\epsilon\neq x)\leq\lambda(x_\epsilon\neq x)+\epsilon\nu^{1/2}(B_i)\leq3\epsilon\nu^{1/2}(B_i).
$$
It follows that if $\epsilon$ is small enough,
$$
\begin{array}{c}
\|Sx\|_X\geq\|Tx_\epsilon\|-\|S(x-x_\epsilon)\|\geq\vspace{0.3cm}\\
\geq\delta\lambda(supp\ x_\epsilon)-\|S\|\|x-x_\epsilon\|_{L_q(\nu)}\geq\vspace{0.3cm}\\
\geq\delta\cdot\lambda(B_i)-\delta\cdot2\epsilon\nu(B_i)^{1/2}
-\|S\|(\frac{3\epsilon}{b_i})^{1/q}\nu(B_i)^{1/2q}\geq\vspace{0.3cm}\\
\geq\delta b_i\cdot\nu(B_i)-3\delta\cdot\epsilon\nu(B_i)^{1/2}
-\|S\|(\frac{3\epsilon}{b_i})^{1/q}\nu(B_i)^{1/2q}>\vspace{0.3cm}\\
>\frac 1 2 \delta b_i\nu(B_i)>0.
\end{array}
$$

So $S:L_q(B_i,\nu)\rightarrow X$ is $l_2$-strictly singular operator but not narrow.

Following the construction in Proposition 3.1 ~\cite{MPRS} for each $\epsilon>0$ we can find tree $\{A_{m,k}\}$ of $B_i$ and operator $\widetilde{S}:L_q(B_i,\Sigma_1,\nu)\rightarrow X,$ where $\Sigma_1$ is a $\Sigma$-algebra generated by $\{A_{m,k}\},$ with the following properties:

$\ (P1) \ \nu(A_{m,k})=2^{-m}\nu(B_i), \ \forall m,k;$

$\ (P2) \ \|\widetilde{S}x\|\geq\frac 1 4 \delta \lambda(B_i)$ for each mean zero sign $x\in L_q(B_i,\Sigma_1,\nu);$

$\ (P3) \ \widetilde{S}(h'_{1})=0$ and $\widetilde{S}h'_{n}=(P_{s_n}-P_{s_{n-1}})Sh'_n,$ where $0=s_1<s_2,\dots,$ $P_n$ are the basis projections in $X,$ and $\{h'_{2^m+k}\}$ is Haar system w.r.t. the tree $\{A_{m,k}\},$ normalized in $L_q(B_i,\Sigma_1,\nu);$

$\ (P4) \ $for all $x\in L_q(\nu), \ \|x\|=1,$ we have $\|\widetilde{S}x\|\leq\|Sx\|+\epsilon;$

$\ (P5) \ $for some sequence $\epsilon_L\rightarrow0$ as $L\rightarrow0$ and each $x\in L_q(\nu), \ \|x\|=1,$ of the form $x=\sum\limits_{n=L}^N\beta_n h'_n$ we have $\|\widetilde{S}x\|\leq\|Sx\|+\epsilon_L.$

Note that if we consider a map $J:\Sigma_1(B_i)\rightarrow \mathfrak{B}(0,1),$ such that
$$
\begin{array}{l}
J(A_{m.k})=\Delta_{m,k}=[\frac{k-1}{2^m},\frac{k}{2^m}]; \\
m(J(A_{m.k}))=\nu^{-1}(B_i)\cdot\nu(A_{m,k}),
\end{array}
$$
where $\mathfrak{B}(0,1)$ is Borel algebra, $m$ is Lebesgue measure, then $L_q(B_i,\Sigma_1,\nu)$ is isometric to $L_q=L_q((0,1),\mathfrak{B},m).$ Thus we have operator $\tilde{S}:L_q\rightarrow X$ of special structure, not narrow and $l_2$-strictly singular. It contradicts with Theorem B ~\cite{MPRS}. $\Box$

It remains an open question whether any Banach space which does not contain isomorphic copies of $l_2$ is Lyapunov.

\end{document}